\documentclass[twoside,leqno]{article}

\usepackage{amssymb}
\usepackage{amsmath}
\usepackage[cp1250]{inputenc}
\usepackage{enumerate}

\newcommand{\punkt}{\hspace{-0.3cm} . \hspace{0.1cm}}

\numberwithin{equation}{section}
\newtheorem{thm}[equation]{\sc Theorem}
\newtheorem{lem}[equation]{\sc Lemma}
\newtheorem{cor}[equation]{\sc Corollary}
\newtheorem{prop}[equation]{\sc Proposition}
\newtheorem{rem}[equation]{\sc Remark }

\makeatletter
\renewcommand{\@seccntformat }[1]{\csname the#1\endcsname. }
\makeatother

 1

\def\NN{\mathbb{N}}

\def\MM{\mathbb{M}}

\def\ZZ{\mathbb{Z}}

\def\QQ{\mathbb{Q}}

\def\CB{{\cal B}}

\def\CE{{\cal E}}
\def\CH{{\cal H}}

\def\CJ{{\cal J}}
\def\CK{{\cal K}}
\def\CM{{\cal M}}

\def\CS{{\cal S}}

\def\rep{\mbox{\rm rep}}
\def\brep{\mbox{\bf rep}}
\def\repf{\mbox{\rm repf}}

\def\vsp{\vspace*{1.5ex}}

\def\epv {{$\mbox{}$\hfill ${\Box}$\vspace*{1.5ex} }}

\def\valg#1{{\bf Alg}_K(d)}

\def\ov{\overline}

\def\bdim{\mbox{\bf dim}}
\def\dimv{\mbox{\bf dim}}
\def\ind{\mbox{{\rm ind}}}
\def\mod{\mbox{{\rm mod}}}

\def\latt{\mbox{{\rm latt}}}

\def\mapr#1#2{\smash{\mathop{\longrightarrow}\limits^{#1}_{#2}}}

\def\mapl#1#2{\smash{\mathop{\longleftarrow}\limits^{#1}_{#2}}}

\def\Hom{\mbox{\rm Hom}}

\def\Ext{\mbox{\rm Ext}}

\def\rad{\mbox{\rm rad}}
\def\soc{\mbox{\rm soc\,}}

\def\K#1{#1^{(K)}}

\def\wt#1{\widetilde{#1}}

\def\K#1{{#1}^{(K)}}

\def\L#1{{#1}^{(L)}}

\def\bc{\begin{center}}
\def\ec{\end{center}}

\def\vsp{\vspace*{1.5ex}}

\begin{document}

\renewcommand{\thefootnote}{}
\def\refname{\begin{center}\normalsize{\it REFERENCES}\end{center}}

 \title{\large{\bf On Lie algebras associated with representation directed algebras}
   \\
   \vspace{0.5em}}
\author{  \normalsize {\bf Stanis{\l}aw Kasjan and Justyna  Kosakowska }
\\
 \small{\it Faculty of Mathematics and Computer Science,
Nicolaus Copernicus University}\\ \small{\it ul. Chopina 12/18,
87-100 Toru\'n, Poland} \\
\small{\it E-mail: skasjan@mat.uni.torun.pl\quad
justus@mat.uni.torun.pl}}
\date{}
\maketitle


\footnote{The second named author is partially supported by Research Grant No. N N201 269 135 of Polish Ministry of Science and High Education}

\begin{abstract}
 Let $B$ be a representation-finite
$\mathbb{C}$-algebra.
 The $\mathbb{Z}$-Lie algebra $L(B)$ associated with $B$ has been defined by Ch. Riedtmann in \cite{riedtmann}.
 If $B$ is representation-directed there is another $\mathbb{Z}$-Lie
 algebra associated with $B$ defined by C. M. Ringel in \cite{ringel1} and denoted by
 $\CK(B)$.

We prove that the Lie algebras $L(B)$ and $\CK(B)$ are isomorphic
for any representation-directed $\mathbb{C}$-algebra $B$.  
\end{abstract}



\normalsize

\section{Introduction}\label{sec:introduction} Let $\mathbb{C}$ be the~field of
complex numbers and let $B$ be a~finite dimensional associative
basic $\mathbb{C}$-algebra with unit element. By a~$B$-module we
mean a~finite dimensional right $B$-module. Assume that $B$ is
a~representation-directed algebra (see Section
\ref{sec:rep-directed} for definitions). In particular, the
algebra $B$ is  representation-finite, that is, there is only
finitely many isomorphism classes of indecomposable $B$-modules.
Let $\CM(B)$ (resp. $\ind(B)$) be a set of representatives of all
isomorphism classes of $B$-modules (resp. indecomposable
$B$-modules). In \cite{riedtmann}, with any~representation-finite
algebra $B$, Ch. Riedtmann associated the $\mathbb{Z}$-Lie algebra
$L(B)$, which is the free $\mathbb{Z}$-module with basis $\{v_X\;
;\; X\in\ind(B)\}$. If $B$ is representation-directed, then
the~Lie algebra structure on $L(B)$ is defined in the following
way. Let $X$, $Y$ be non-isomorphic indecomposable $B$-modules
such that $\Ext_B^1(X,Y)=0$. Put
$$\CE(X,Y;Z)=\{\; U\subseteq Z\; ; \; \; U\cong X,\; Z/U\cong Y\}.$$
This is locally a closed subset of a~product of Grassmann
varieties, see \cite{riedtmann} and Section \ref{sec:riedmann} for
details. We set
$$
[v_X,v_Y]=\left\{\begin{array}{ll} \chi(\CE(X,Y;Z))\cdot v_Z &
\mbox{if there is an~indecomposable } B\mbox{-module } Z\\
&\mbox{and  a~short exact sequence } \\ & 0\to X\to Z\to Y\to 0,
\\ & \\ 0 & \mbox{otherwise,}
\end{array} , \right.
$$
where $\chi(\CE(X,Y;Z))$ denotes the Euler-Poincar\'{e}
characteristic of $\CE(X,Y;Z)$.  This formula defines the unique
Lie algebra structure on $L(B)$.

On the other hand, in \cite{ringel1}, C. M. Ringel  associated
with a~representation-directed algebra $B$ the $\mathbb{Z}$-Lie
algebra $\CK(B)$, which is the free $\mathbb{Z}$-module with basis
$\{u_X\; ;\; X\in\ind(B)\}$. If $X$, $Y$ are non-isomorphic
indecomposable $B$-modules such that $\Ext_B^1(X,Y)=0$, we set
$$
[u_Y,u_X]=\left\{\begin{array}{ll} \varphi_{YX}^Z(1)\cdot u_Z &
\mbox{if there is an~indecomposable } B\mbox{-module } Z\\
&\mbox{and a~short exact sequence } \\ & 0\to X\to Z\to Y\to 0
\\ & \\ 0 & \mbox{otherwise,}
\end{array} , \right.
 $$  where $\varphi_{YX}^Z$
are Hall polynomials (see \cite{ringel1} and Section
\ref{sec:ringel} for details). This formula defines the unique Lie
algebra structure on  $\CK(B)$.

If $B=\mathbb{C}Q$ is the~path algebra of a~Dynkin quiver $Q$,
then $\mathbb{C}\otimes_{\mathbb{Z}}L(B)\cong {\bf n}_+(Q)$ and
$\mathbb{C}\otimes_{\mathbb{Z}}\CK(B)\cong {\bf n}_+(Q)$, where
${\bf n}_+(Q)$ is the positive part of the complex semisimple Lie
algebra of type $Q$ (see \cite[Proposition 5.3]{riedtmann} and
\cite[Corollary 3]{ringel90}). Therefore for path algebras
$B=\mathbb{C}Q$ of Dynkin quivers the Lie algebras $L(B)$ and
$\CK(B)$ are isomorphic.  Ch. Riedtmann conjectured  in
\cite{riedtmann} that the Lie algebras $L(B)$ and $\CK(B)$ are
isomorphic for any representation-directed $\mathbb{C}$-algebra
$B$. The main aim of this paper is to prove this conjecture, more
precisely we prove the following theorem.

\begin{thm}\punkt Let $B$ be a~representation-directed $\mathbb{C}$-algebra.
There is a~Lie algebra isomorphism
$$
F:L(B)\mapr{}{}\CK(B)
$$
given by the formula
$$
F(u_{X})=(-1)^{(\dim_{\mathbb{C}}X- 1)} v_{X}.
$$\label{thm:main}\end{thm}

Theorem \ref{thm:main} is an immediate consequence of the
following fact.

\begin{thm}\punkt Let $B$ be a~representation-directed $\mathbb{C}$-algebra
and let $X$, $Y$, $Z$ be $B$-modules {\rm (}not necessarily indecomposable{\rm )}. Then
$$ \chi(\CE(X,Y;Z))=\varphi_{YX}^Z(1).$$ \label{thm:main1}\end{thm}

The reason for the appearance of $(-1)^{(\dim_{\mathbb{C}}X- 1)}$
in the formula in \ref{thm:main} is that the original definitions
of the lie brackets in $L(B)$ and $\CK(B)$ are not "compatible" in
the sense that $[v_X,v_Y]$ refers to the extensions of $Y$ by $X$,
whereas $[u_X,u_Y]$ to the extensions of $X$ by $Y$.


The paper is organized as follows.

\begin{itemize}
\item In Section \ref{sec:rep-directed}
we collect some facts of the representation theory of algebras. In
particular, we recall basic properties of representation-directed
algebras and their Auslander-Reiten quivers.
\item In Section \ref{sec:riedmann} we recall definitions and
basic properties of the $\mathbb{Z}$-Lie algebra $L(B)$.
\item In Section \ref{sec:ringel} we recall definitions and
basic properties of the $\mathbb{Z}$-Lie algebra $\CK(B)$, Hall
polynomials and Ringel-Hall algebras.
\item In Section \ref{sec:Z-form}  we show that the modules over a
representation-directed algebra can be "defined over $\ZZ$"
independently on the base field.
\item In Section \ref{sec:grassmann} we investigate the subset
$\CE(X,Y;Z)$ of a~product of Grassmann varieties. Results of Section \ref{sec:Z-form} and \ref{sec:grassmann} allow
to describe $ \chi(\CE(X,Y;Z))$ in terms of Hall polynomials $\varphi_{YX}^Z.$
\item Section \ref{sec:proof} contains proofs of Theorems \ref{thm:main} and \ref{thm:main1}
and consequences of the theorems.
\end{itemize}

The motivation for the study of Lie algebras $L(B)$ and $\CK(B)$
is their connection with Ringel-Hall algebras, generic extensions,
and quantum groups (see \cite{reineke2001}, \cite{ringel},
\cite{ringel1}, \cite{ringel89} and \cite{ringel90}). Hall
polynomials and the Euler-Poincar\'{e} characteristic $\chi$
are developed in \cite{CC}, where the authors
investigate connections between indecomposable representation of
Dynkin quivers of type $\mathbb{A}$, $\mathbb{D}$, $\mathbb{E}$
and the cluster variables of cluster algebras of type
$\mathbb{A}$, $\mathbb{D}$, $\mathbb{E}$.

Results of this paper where presented by the second named author
during the conference "Colloque d'algebre non-commutative" in University
of Sherbrooke (June 2008) and on the "Seminar Darstellungstheorie" in Bielefeld
during her stay in University of Bielefeld supported by SFB (September 2008).

\section{Preliminaries}\label{sec:rep-directed}

In this section we shortly recall basic definitions, notation and
facts of the representation theory of representation-directed
algebras. For the basic concepts of representation theory the
reader is referred to \cite{ASS} and \cite{ARS}.

We consider algebras of the form $B=KQ/I$, where $Q=(Q_0,Q_1)$ is
a finite quiver and $I$ is an admissible ideal of the path algebra
$KQ$ of $Q$. Recall, that an ideal $I$ is {\bf admissible} if it
is contained in the ideal generated by the paths of length at
least 2 and there is a number $N$ such that every path of length
$N$ belongs to $I$. If $Q$ is {\bf acyclic}, that is, contains no
oriented cycles, the latter condition is satisfied automatically.

If $B$ is as above then it is a finite dimensional associative
basic $K$-algebra. All modules  are assumed to be right finite
dimensional. Denote by $\mod(B)$ the category of all such
$B$-modules. The category $\mod(B)$ is equivalent to the category
$\rep_K(Q,I)$ of finite dimensional $K$-representations of $Q$
satisfying the relations from $I$, see \cite[Definition 1.4]{ASS}.
We identify the two categories.

Given a field $K$ and a dimension vector ${\bf d}\in \NN^{Q_0}$
let $\brep_K(Q,{\bf d})$ be the variety of $K$-representations of
$Q$ with the dimension vector ${\bf d}$. That is,
$$
\brep_K(Q,{\bf d})=\prod_{\alpha\in Q_1}\MM_{{\bf
d}_{t(\alpha)},{\bf d}_{s(\alpha)}}(K),
$$
where $t(\alpha)$ and $s(\alpha)$ denote the terminus and the
source of $\alpha$, respectively, and $\MM_{a,b}(K)$ is the space
of all $a\times b$-matrices. We identify a point of
$\brep_K(Q,{\bf d})$ with the corresponding representation of $Q$
in the usual way.

The group $Gl({\bf d}, K)=\prod_{x\in Q_0}Gl({\bf d}_x, K)$ acts
"by conjugations" on $\brep_K(Q,{\bf d})$ and the orbit of the
point $M$ is the set of points corresponding to the
representations isomorphic to $M$.

For any $B$-module $M$, identified with the representation
$(M_x,M_{\alpha})_{x\in Q_0,\alpha\in Q_1}$, denote by $\bdim\,
M\in \mathbb{N}^{Q_0}$ \textbf{the dimension vector} of $M$
defined by   $(\bdim\, M)_x=\dim_KM_x$, for any $x\in Q_0$.

Let $\CM(B)$ (resp. $\ind(B)$) be a set of representatives of all
isomorphism classes of $B$-module (resp. indecomposable
$B$-modules).

Denote by $\Gamma_B=((\Gamma_B)_0,(\Gamma_B)_1)$ the
Auslander-Reiten quiver of $B$ and  by $\tau_B$ the
Auslander-Reiten translation in $\Gamma_B$. Recall that there is
a~bijection $(\Gamma_B)_0\leftrightarrow\ind(B)$.

Given two indecomposable $B$-modules $X$, $Y$ we write $X\le Y$ if
there is a sequence of nonzero maps
$X=Y_0\mapr{}{}Y_1\mapr{}{}...\mapr{}{}Y_t=Y$ with $Y_0,...,Y_t$
indecomposable. Following \cite{SkWe}, we say that a~$B$-module
$X$ is {\bf directing} if there do not exist indecomposable direct
summands $X_1$, $X_2$ of $X$ and an indecomposable nonprojective
module $Y$ such that $X_1\le \tau_BY$ and $Y\le X_2$.

An~algebra $B$ is said to be {\bf  representation-finite} (resp.
{\bf representation-directed}), if there is only finitely many
isomorphism classes of indecomposable $B$-modules (resp. any
indecomposable $B$-module is directing). It is well known that any
representation-directed algebra is representation-finite
\cite[page 78]{Ri1099}. If $B$ is a connected
representation-directed algebra then the Auslander-Reiten quiver
of $B$ coincides with its (unique) preprojective component.
Moreover, if $B$ is representation-directed then we may enumerate
the indecomposable $B$-modules $X_1,\ldots,X_m$ in such a~way that
$\Ext^1_B(X_i,X_j)=0$, if $i\leq j$, and $\Hom_B(X_i,X_j)=0$, if
$j<i$.

If in addition $B$ is of the form $KQ/I$ and $I$ is an admissible
ideal, then every arrow of $\Gamma_B$ has trivial valuation
\cite[Section VII]{ARS}.

Assuming this, let $P_x^B$ be an indecomposable projective
$B$-module associated with the vertex $x\in Q_0$ and decompose the
radical $\rad\,P^B_x$ of $P^B_x$  into indecomposable direct
summands:
$$
\rad\,P^B_x\cong M_x^1\oplus...\oplus M_x^{r_x}.$$
 We denote the
set $\{\bdim\,M_x^1,...,\bdim\,M_x^{r_x}\}$ of their dimension
vectors by $R^B_x$. Observe that these dimension vectors are
pairwise different. \vsp

\begin{lem}\punkt Let $Q$ be a connected  acyclic quiver, $B=KQ/I$,
$C=LQ/J$, for some fields $K$, $L$ and admissible ideals $I$ and
$J$ of $KQ$ and $LQ$ respectively. Assume that
$\bdim\,P^B_x=\bdim\,P^C_x$, $R^B_x=R^C_x$, for any $x\in Q_0$.

If $B$ is representation-directed then $C$ is also
representation-directed and there is an isomorphism
$$
\sigma:\Gamma_B\mapr{}{}\Gamma_C
$$
of translation quivers such that the dimension vectors of the
$B$-module corresponding to the vertex $m$ of $\Gamma_B$ and the
$C$-module corresponding to the vertex $\sigma(m)$ of $\Gamma_C$
are equal, for any $m$.\label{lem:rep-dir}\end{lem}

{\bf Proof.} The input data of the usual algorithm for
constructing preprojective components of algebra ("knitting
procedure") are the dimension vectors of

(i) the indecomposable projective modules and

(ii) the indecomposable direct summands of the radicals of the
indecomposable projective modules

\noindent (see \cite{DP}, \cite{KaPe}). We mean here an algorithm
determining the "combinatorial data" of preprojective components:
the dimension vectors of indecomposable preprojective modules and
the numbers of arrows between the  vertices.

Applying the algorithm we  construct the Auslader-Reiten quiver of
$B$ (equal to its unique preprojective component). Since the input
data of the algorithm are the same for both algebras, then there
exist unique preprojective component ${\cal P}$ of $C$ and
an~isomorphism $\sigma \Gamma_B\mapr{}{}{\cal P}$. Moreover
$\sigma$ preserves dimension vectors.   The  component ${\cal P}$
is finite thus coincides with $\Gamma_C$ by \cite[Theorem
2.1]{ARS}, and the lemma follows. \epv

In order to investigate modules over representation directed
algebras independently on the base field we need the concept of
lattices over orders.

 By a~$\mathbb{Z}$-{\bf order} we mean a
unitial ring $A$ which is free and finitely generated as a
$\mathbb{Z}$-module. An~{\bf $A$-lattice} is a~right $A$-module,
free and finitely generated as a~$\mathbb{Z}$-module. We denote by
$\latt(A)$ the category of $A$-lattices. Given a~unitial
commutative ring $R$ let $A^{(R)}$ denote the $R$-algebra
$A\otimes_{\mathbb{Z}}R$. We call it {\bf specialization} of $A$
to $R$. Similarly, given an $A$-lattice $X$ we denote by $X^{(R)}$
the $A^{(R)}$-module $A\otimes_{\mathbb{Z}}R$.

We are interested in $\mathbb{Z}$-orders of the form $\mathbb{Z}
Q/I$, where $\mathbb{Z} Q$ is the path $\mathbb{Z}$-algebra of
some quiver $Q$ and $I$ is an admissible ideal of $\mathbb{Z} Q$.
The ring $\mathbb{Z} Q/I$ is a $\mathbb{Z}$-order if an only if
$I$ is a pure $\mathbb{Z}$-submodule of $\mathbb{Z} Q$.

If $A=\mathbb{Z} Q/I$ is a $\mathbb{Z}$-order and $L$ is a field,
then the $L$-algebra $\L{A}$ is isomorphic to $LQ/\L{I}$, where
$I$ is the ideal of $LQ$ generated by the image of $I$ under the
canonical homomorphism $\mathbb{Z} Q\mapr{}{} LQ$. Observe that if
$I$ is admissible, then $\L{I}$ is also.

Given a quiver $Q$ let $\rep_{\mathbb{Z}}(Q)$ be the category of
all finitely generated $\mathbb{Z}$-representations of $Q$, that
is the systems $(N_x,N_{\alpha})_{x\in Q_0,\alpha\in Q_1}$, where
$N_x$ are finitely generated $\mathbb{Z}$-modules and
$N_{\alpha}:N_{s(\alpha)}\mapr{}{}N_{t(\alpha)}$ are
$\mathbb{Z}$-homomorphisms, for $\alpha\in Q_1$. The morphisms are
defined in the usual way.

We denote by $\repf_{\mathbb{Z}}(Q)$ the full subcategory of
$\rep_{\mathbb{Z}}(Q)$ consisting of the representations
$(N_x,N_{\alpha})_{x\in Q_0,\alpha\in Q_1}$ such that $N_x$ is a
free $\mathbb{Z}$-module, for any $x\in Q_0$.

Given  an ideal $I$ of $\mathbb{Z} Q$ we denote by
$\rep_{\mathbb{Z}}(Q,I)$ and $\repf_{\mathbb{Z}}(Q,I)$ the full
subcategories of $\rep_{\mathbb{Z}}(Q)$ and
$\repf_{\mathbb{Z}}(Q)$, respectively, consisting of the
representations satisfying all relations in $I$.

It is clear that if $A=\mathbb{Z} Q/I$ is a $\mathbb{Z}$-order
then there is an equivalence of categories
$$
\latt(A)\cong\repf_{\mathbb{Z}}(Q,I).
$$
>From now on we identify the two categories.

Given an $A$-lattice $N$ corresponding to the representation
$(N_x,N_{\alpha})_{x\in Q_0,\alpha\in Q_1}$ we denote by $\bdim
N\in\NN^{Q_0}$ the dimension vector of $N$, that is $(\bdim N)_x$
equals the rank of $N_x$, $x\in Q_0$.

\section{Lie algebras associated with representation directed algebras}

Given two natural numbers $e\le d$ and a field $K$ we denote by
$Gr_e(d,K)$ the Grassmann variety of $e$-dimensional subspaces of
$K^d$, embedded into the projective space $\mathbb{P}^{{d\choose e}-1}$
via Pl\"ucker embedding.

Let $Q=(Q_0,Q_1)$ be a finite quiver and let $K$ be a~field. Given
two vectors ${\bf e}, {\bf d}\in\NN^{Q_0}$ such that ${\bf e}_x\le
{\bf d}_x$, for any $x\in Q_0$, we set
$$
Gr_{{\bf e}}({\bf d},K)=\prod\limits_{x\in Q_0}Gr_{{\bf e}_x}({\bf
d}_x,K).
 $$

Let $M=(M_x,M_{\alpha})_{x\in Q_0,\alpha\in Q_1}$ be a
$K$-representation of $Q$ with the dimension vector ${\bf d}$. We
identify $M_x$ with  $K^{{\bf d}_x}$, for $x\in Q_0$. Let
$$
Sub_{{\bf e}}(M)\subseteq Gr_{{\bf e}}({\bf d},K)
$$
be the set of the tuples $N=(N_x)_{x\in Q_0}$ of subspaces which
form are subrepresentation of $M$, that is, such that
$M_{\alpha}(N_{s(\alpha)})\subseteq N_{t(\alpha)}$, for any
$\alpha\in Q_1$. Given $N\in Sub_{{\bf e}}(M)$ we treat $N$ as a
representation of $Q$ in the natural way.

Assume that $i_x$ indicates a homogeneous coordinate in
$\mathbb{P}^{{{\bf d}_x\choose {\bf e}_x}-1}$, for $x\in Q_0$, and
put $\underline{i}=(i_x)_{x\in Q_0}$. We denote by ${\cal
U}_{\underline{i}}$ the open subset of $Gr_{{\bf e}}({\bf d},K)$
consisting of the tuples $N=(N_x)_{x\in Q_0}$ of subspaces such
that the $i_x$th coordinate of $N_x$ is nonzero, for $x\in Q_0$.

The following assertion follows by standard arguments. \vsp

\begin{lem}\punkt We keep the notaion introduced above.

{\rm (1)} The set $Sub_{{\bf e}}(M)$ is a Zariski-closed subset of
$Gr_{{\bf e}}({\bf d},K)$.

{\rm (2)} For any possible $\underline{i}$, there exist regular
maps
$$
\begin{array}{rcccl}
\brep_K(Q,{\bf e})&\mapl{\eta_{\underline{i}}}{}&Sub_{{\bf
e}}(M)\cap {\cal U}_{\underline{i}}&\mapr{\rho_{\underline{i}}}{}&
\brep_K(Q,{\bf d}-{\bf e})
\end{array}
$$
such that a subrepresentation $N$ of $M$ corresponds to
$\eta_{\underline{i}}(N)$ and the factor-representation $M/N$
corresponds to $\rho_{\underline{i}}(N)$, for any $N\in Sub_{{\bf
e}}(M)\cap {\cal U}_{\underline{i}}$.
\label{lem:grassmann}\end{lem}

 Let $B=KQ/I$ for some admissible ideal $I$ of $KQ$, ${\bf e}\in \NN^{Q_0}$, and let
$N_1,N_2,M$ be $B$-modules identified with representations of $Q$.
We set
\begin{equation} \CE(N_1,N_2;M)=\{\; U\subseteq M\; ;
\; U\in \mod(B) ,\; U\cong N_1,\; X/U\cong N_2\}. \label{eq:V}
\end{equation}
Note that $\CE(N_1,N_2;M)$ is empty unless $\bdim\, M=\bdim\,N_1+\bdim\,N_2$.
The set $\CE(N_1,N_2;M)$ is a subset of
$$
Gr_{\scriptsize {\bdim\,N_1}}(\bdim\,M,K).
$$

\begin{lem}\punkt Let $K$ be an~algebraically closed field. Then
$\CE(N_1,N_2;M)$ is a~locally closed subset of
$Gr_{\scriptsize {\bdim\,N_1}}(\bdim\,M,K)$.
\label{lem:locclosed}\end{lem}

{\bf Proof.} The assertion follows from the well-known fact that
an~orbit of an algebraic group action is locally closed by Lemma
\ref{lem:locclosed}.\vsp \epv

\subsection{Riedtmann's Lie algebras}\label{sec:riedmann}

Let $B=\mathbb{C}K/I$ be a representation-finite
$\mathbb{C}$-algebra. In \cite{riedtmann}, Ch. Riedtmann
associated with  $B$ the Lie algebra $L(B)$ as follows. Let
$\chi(\CE(N_1,N_2;M))$ denote the Euler-Poincar\'{e}
characteristic of a~locally closed subset $\CE(N_1,N_2;M)$ defined
in \ref{eq:V}. Let $R(B)$ be the~free $\mathbb{Z}$-module with
basis $\{ v_M\; ;\; M\in \CM(B) \}$.  The formula
\begin{equation} v_M\cdot v_N=\sum_{X\in
\CM(B)}\chi(\CE(M,N;X))v_X\label{eq:Riedtmann-alg}\end{equation}
defines an associative $\ZZ$-algebra structure on $R(B)$
(see \cite[2.3]{riedtmann}). Note that the sum in (\ref{eq:Riedtmann-alg}) is finite, since
there is only finitely many isomorphism classes of modules $X$
with $\bdim\, X= \bdim\,N+\bdim\,M$ and $\chi(\CE(M,N;X))=0$ if
$\bdim\, X\neq \bdim\,N+\bdim\,M$. The $\mathbb{Z}$-submodule
\begin{equation} L(B)=\bigoplus_{X\in
\begin{scriptsize}\ind\end{scriptsize}(B)}\ZZ v_X\label{eq:L}\end{equation} of $R(B)$ is a~Lie
subalgebra of finite rank $|\ind(B)|$ of the Lie algebra $R(B)$
with respect to the Lie bracket $[x,y]=xy-yx$, where $|S|$ denotes
the cardinality of a~finite set $S$ (see \cite[2.3]{riedtmann}).

\subsection{Ringel's Lie algebras}\label{sec:ringel}

Following ideas of C. M. Ringel from \cite{ringel1}, we define
Hall and Lie algebras of representation-directed algebras. Assume
that $\Gamma=\Gamma_B$  is the  Auslander-Reiten quiver of a
representation-directed $K$-algebra $B=KQ/I$.  For any $x\in
\Gamma_0$, we denote by $M(B,x)\in\ind(B)$ the~indecomposable
$B$-module corresponding to $x$. Denote by $\CB(\Gamma)$ the set
of all functions $a:\Gamma_0\to \mathbb{N}$. With any~function
$a:\Gamma_0\to\mathbb{N}$, we associate the $B$-module
$M(B,a)=\bigoplus_{x\in\Gamma_0}M(B,x)^{a(x)}$. This establishes
a~bijection between the set $\CB(\Gamma)$ and the set $\CM(B)$ of
all isomorphism classes of $B$-modules.

Let $K$ be a~finite field, $C$ be a~$K$-algebra and $N_1,N_2,M\in
\mod(C)$. We set
\begin{equation}
F_{N_2,N_1}^M=|\CE(N_1,N_2;M)|,\label{eq:F}\end{equation}
where $|X|$ denotes the cardinality of a~finite set $X$.

The following theorem is due to Ringel (\cite[Theorem
1]{ringel1}).

\begin{thm}\punkt Let $\Gamma$ be a~directed Auslander-Reiten quiver $($i.e.
there exists a~representation directed $\mathbb{C}$-algebra $B$
with $\Gamma_B=\Gamma)$, and $a,b,c\in \CB(\Gamma)$. There exists
a~polynomial $\varphi_{ca}^b\in\mathbb{Z}[T]$ $($called {\bf Hall
polynomial}$)$ with the following property: if $K$ is a~finite
field, and $C$ is a~$K$-algebra with $\Gamma=\Gamma_C$, then
$$F_{M(C,c),M(C,a)}^{M(C,b)}=\varphi_{ca}^b(|K|). $$
\label{thm:hall-poly}\end{thm}

Let $\Gamma$ be a~directed Auslander-Reiten quiver. Following
\cite{ringel1}, we define the {\bf degenerate Ringel-Hall algebra}
$\CH(\Gamma)_1$ to be the~free $\mathbb{Z}$-module with basis
$\{u_a\}_{a\in \CB(\Gamma)}$ and multiplication given by the
formula
\begin{equation} u_cu_a=\sum_{b\in
\CB(\Gamma)}\varphi_{ca}^b(1)u_b.\label{eq:ringel-alg}\end{equation}
By \cite[Proposition 5]{ringel1}, $\CH(\Gamma)_1$ is
an~associative $\mathbb{Z}$-algebra with unit element and the
$\mathbb{Z}$-submodule \begin{equation}
\CK(\Gamma)=\bigoplus_{x\in \Gamma_0}\ZZ
u_x\label{eq:K}\end{equation} of $\CH(\Gamma)_1$ is a~Lie
subalgebra of finite rank $|\Gamma_0|=|\ind(B)|$ of the Lie
algebra $\CH(\Gamma)_1$ (with the Lie bracket $[x,y]=xy-yx$).

Assume that $B$ be a~representation-directed algebra over
the~field $\mathbb{C}$ of complex numbers. Let $\Gamma_B$ be the
Auslander-Reiten quiver of $B$. We set $\CH(B)_1:=\CH(\Gamma_B)_1$
and $\CK(B):=\CK(\Gamma_B)$.

\section{ Lattices} \label{sec:Z-form}

In this section we prove that for any representation-directed
algebra $B=KQ/I$ there is a $\mathbb{Z}$-order $A$ such that the
$K$-algebras  $\K{A}$ and $B$ are isomorphic and the
specialization $\L{A}$ of $A$ to any other field $L$ is
representation-directed. Moreover, the $\L{A}$-modules are
induced from some $A$-lattices chosen independently on $L$.

Recall, that a representation $(N_x,N_{\alpha})_{x\in
Q_0,\alpha\in Q_1}$of $Q$ over a field $L$ is called thin if
$\dim_LN_x\le 1$, for $x\in Q_0$. We say that a
$\mathbb{Z}$-representation $(N_x,N_{\alpha})_{x\in Q_0,\alpha\in
Q_1}$ of $Q$ is {\bf thin} if $N_x\cong \mathbb{Z}$ or $N_x=0$,
for $x\in Q_0$ and $N_{\alpha}$ is either 0 or invertible, for any
arrow $\alpha$.

If $A=\ZZ Q/I$ is an order and $I$ is an admissible ideal then we
say that an $A$-lattice $N$ is {\bf thin} provided so is the
corresponding representation of $Q$.

\vsp

\begin{lem}\punkt  Let $A=\mathbb{Z} Q/I$ be a $\mathbb{Z}$-order, where $I$ is
an admissible ideal of $\mathbb{Z} Q$ and assume that $\L{A}$ is a
representation-directed algebra, for any field $L$. Assume that
$X$ is an $A$-lattice such that $\L{X}$ is an indecomposable thin
$\L{A}$-module, for any field $L$, and $Y$ is a thin
indecomposable lattice. If $\dimv(\K{X})=\dimv (\K{Y})$ for some
field $K$, then $X\cong Y$ as $A$-lattices. \label{lem:thin}
\end{lem} \vsp

{\bf Proof.} Given a thin representation $Z=(Z_x,Z_{\alpha})_{x\in
Q_0,\alpha\in Q_1}$ of $Q$  over a field let $N(Z)=(N_0,N_1)$ be
the subquiver of $Q$ with the set of vertices $N_0=\{x\in
Q_0:Z_x\neq 0\}$ and arrows $N_1=\{\alpha\in Q_1:Z_{\alpha}\neq
0\}$. Note that $Z$ is indecomposable if and only if $N(Z)$ is
connected. It follows that the $\L{A}$-module $\L{Y}$ is
indecomposable, for any field $L$. Therefore $\L{X}\cong \L{Y}$,
since the dimension vector determines uniquely the isomorphism
class of an indecomposable module over a representation-directed
algebra \cite{bo83}, \cite[Chap. IX, Prop. 3.1]{ASS}. The quiver
$N(\L{X})$ does not depend on the field $L$, hence
 $X$ is a thin lattice. Now it is easy to observe that the lattices
$X$ and $Y$ are isomorphic (apply the isomorphism $\L{X}\cong
\L{Y}$ for a field $L$ of characteristic 0). \epv

The following theorem is a version of Lemma 4.3 from
\cite{Ka}.\vsp

\begin{thm}\punkt Let $B$ be a representation-directed algebra
of the form $KQ/J$ for an acyclic quiver $Q$ and an admissible
ideal $J$ in $KQ$. There exist:

{\rm (a)} a $\mathbb{Z}$-order $A=\mathbb{Z} Q/I$, where $I$ is a
two-sided admissible ideal $I$ of $\mathbb{Z} Q$.

{\rm (b)} $A$-lattices $X_1,...,X_r$,

{\rm (c)} a set $\CJ\subseteq \{1,...,r\}^2$ and a family of
$A$-homomorphisms $f_{ij}:X_i\mapr{}{} X_j$, for $(i,j)\in\CJ$,

\noindent such that

{\rm (1)} $A^{(K)}\cong B$,

\noindent
and, for any field $L$:

{\rm (2)} $A^{(L)}$ is a representation-directed $L$-algebra,

{\rm (3)} $X_1^{(L)},...,X_r^{(L)}$ are indecomposable pairwise
nonisomorphic $A^{(L)}$-modules,

{\rm (4)} $f_{ij}^{(L)}:X_i^{(L)}\mapr{}{}X_j^{(L)}$ are
irreducible maps, for $(i,j)\in\CJ$, and if there is an
irreducible map $X_i^{(L)}\mapr{}{}X_j^{(L)}$ for some $i,j$ then
$(i,j)\in\CJ$.

{\rm (5)} the Auslander-Reiten quiver $\Gamma_{A^{(L)}}$ of
$A^{(L)}$ has vertices $$x_1=[X_1^{(L)}],...,x_r=[X_r^{(L)}]$$ and
there is an arrow $x_i\mapr{}{}x_j$ if and only if $(i,j)\in\CJ$,

{\rm (6)} the dimension vector ${\bf dim}(X_i^{(L)})$ does not
depend on $L$, for any fixed $i$. \label{thm:Z-form} \end{thm}\vsp

{\bf Proof.} The algebra $B$ is schurian. Thanks to multiplicative
basis theorem (it is enough to apply the easier "triangular" case,
see \cite{Bo1}) we may assume that $J$ is generated by
zero-relations and commutativity relations. Let $\mathfrak{Z}$ be
the set of all paths in $Q$ belonging to $J$ and let
$\mathfrak{C}$ denote the set of all elements $u-w\in J$, where
the paths $u$, $w$ do not belong to $\mathfrak{Z}$. Clearly, $J$
is the $K$-subspace of $KQ$ spanned by $\mathfrak{Z}\cup
\mathfrak{C}$. For any pair $x,y\in Q_0$ such that $e_xBe_y\neq 0$
choose a path $u_{x,y}\notin J$ starting at $x$ and ending at $y$.
Let $\mathfrak{B}$ be a set of all the paths $u_{x,y}$.  Then the
$J$-cosets of the elements of $\mathfrak{B}$ form a $K$-basis of
$B$.

Let $I$ be the $\mathbb{Z}$-submodule of $\mathbb{Z} Q$ spanned by
$\mathfrak{Z}\cup \mathfrak{C}$. It is easy to check that $I$ is a
two-sided ideal of $\mathbb{Z} Q$. Moreover, $A=\mathbb{Z} Q/I$ is
isomorphic to
$$
\bigoplus_{u\in\mathfrak{B}}\ZZ u
$$
as a $\ZZ$-module, therefore $A$ is a $\ZZ$-order.

It is clear that $A^{(K)}\cong B$.

We say that a set $\CS=\{Z_1,...,Z_l\}$ of pairwise nonisomorphic
$B$-modules is {\bf properly defined over } $\mathbb{Z}$ if there exist
$A$-lattices $\wt{Z_1},...,\wt{Z_l}$, a set $\CJ_{\CS}\subseteq
\{1,...,l\}^2$ and a family of $A$-homomorphisms
$g_{ij}:\wt{Z_i}\mapr{}{} \wt{Z_j}$, for $(i,j)\in\CJ_{\CS}$,
 such that  $\K{\wt{Z_i}}\cong Z_i$, for $i=1,...,l$, and the conditions
(3), (4), (6) of the theorem are
 satisfied, for any field $L$,
 with $X_1,...,X_r$, $f_{ij}$, $(i,j)\in\CJ$  interchanged by
 $\wt{Z_1},...,\wt{Z_l}$, $g_{ij}$, $(i,j)\in\CJ_{\CS}$,
 respectively. In this case we call $\{\wt{Z_1},...,\wt{Z_l}\}$ a
 {\bf set of $\mathbb{Z}$-frames of} $\CS$.

We need to prove that $\ind(B)$ is properly defined over
$\mathbb{Z}$.

Let $P_x=P^B_x$ be the indecomposable projective $B$-module
associated with the vertex $x$ of $Q$ and assume that $\rad
P_x\cong M^1_{x}\oplus...\oplus M_{x}^{r_x}$ is a decomposition of
$\rad P_x$ into indecomposable direct summands. It is easy to see
that the set $\{P_x,M_{x}^{i}\; ; \;x\in Q_0, i=1,...,r_x\}$ is
properly defined over $\mathbb{Z}$ and there is a set
$\{\wt{P}_x,\wt{M}_{x}^{i}\; ; \;x\in Q_0, i=1,...,r_x\}$ of
$\mathbb{Z}$-frames consisting of thin lattices. We assume that
the notation is the natural one, that is, $\wt{P}^{(K)}_x\cong
P_x$, etc. Then, for any field $L$, $\wt{P}^{(L)}_x$ is the
indecomposable projective $\L{A}$-module associated with the
vertex $x$, $\L{(M_{x}^{i})}$ are indecomposable $\L{A}$-modules
and
$$
\L{P}_x\cong \L{(M_{x}^{1})}\oplus...\oplus \L{(M_{x}^{r_x})}.
$$

Similarly, the set  of representatives $I_x$, $x\in Q_0$, of the
isomorphism classes of indecomposable injective $B$-modules
together with indecomposable direct summands of $I_x/\soc I_x$ is
properly defined over $\mathbb{Z}$.

It follows by Lemma \ref{lem:rep-dir} that the algebra $\L{A}$ is
representation-directed, for any field $L$. Moreover, the
combinatorial data of the Auslander-Reiten quiver of $\L{A}$ do
not depend on the field $L$. It remains to prove that the vertices
of those quivers can be realized by lattices over $A$.

Inductively we construct  sets  $\CS_i\subseteq \ind(B)$, such that

(i) $\CS_0=\{P\}$, where $P$ is a simple projective $B$-module,

(ii) $\CS_l=\ind(B)$, for some natural number $l$,

\noindent and, for any $i=0,...,l-1$:

(iii) $\CS_i\subseteq \CS_{i+1}$,

(iv) $\CS_i$ is properly defined over $\mathbb{Z}$,

(v) $\CS_i$ is closed under predecessors in $\Gamma_B$.

Suppose that the set $\CS_n$ is already defined for some $n$ and
$\CS_n\neq \ind(B)$. Let $Z\in\ind(B)\setminus \CS_n$ be such that
every proper predecessor of $Z$ belongs to $\CS_n$.

If $Z$ is projective, say $Z=P_x$, then
$M_{x}^{1},...,M_{x}^{r_x}\in\CS_n$ up to isomorphism. Thus there
are lattices $N_1,...,N_{r_x}$ in the set $\wt{\CS}_n$ of
$\mathbb{Z}$-frames of $\CS_n$ such that $\K{N}_j\cong M_{x}^{j}$,
for any $j=1,...,r_x$. By Lemma \ref{lem:thin}, we conclude that
$N_j\cong \wt{M}_{x}^{j}$, for $j=1,...,r_x$.
 We
set $\CS_{n+1}={\CS}_n\cup\{{Z}\}$ and we see that
$\wt{\CS}_n\cup\{\wt{P_x}\}$ is a set of $\mathbb{Z}$-frames of
$\CS_{n+1}$.

If $Z$ is not projective, then  let
$$
0\mapr{}{}X\mapr{\eta}{}\bigoplus_{i=1}^{m_Z}E_i\mapr{\nu}{}Z\mapr{}{}0
$$
 be the Auslander-Reiten sequence terminating at $Z$,
 where $E_i$ are indecomposable modules, for $i=1,...,m_Z$.
 We denote by $\eta_i$ and $\nu_i$ the component maps of $\eta$ and $\nu$,
 respectively, for $i=1,...,m_Z$.

Let $\wt{X}$ and $\wt{E_i}$ be $A$-lattices in $\wt{\CS_n}$ such
that $\K{\wt{X}}\cong X$ and $\K{\wt{E_i}}\cong E_i$, for all $i$.
There are nonzero $A$-homomorphisms
$\wt{\eta_i}:\wt{X}\mapr{}{}\wt{E_i}$ such that the induced map
$\L{\wt{\eta_i}}:\L{\wt{X}}\mapr{}{}\L{\wt{E_i}}$ is irreducible
for every $i$. Since  every arrow in $\Gamma_B$ (and consequently
in $\Gamma_{\L{A}}$, for any field $L$) has trivial valuation,
then
$$
\L{X}\mapr{\L{\wt{\eta}}}{}\bigoplus_{i=1}^{m_Z}\L{E_i}, \leqno(*)
$$
where $\wt{\eta}=[\wt{\eta}_i]_{i=1,...,m_Z}$, is a left minimal
almost split map, for any field $L$.

Moreover, $\L{X}$ is not injective, hence the  map $\L{\wt{\eta}}$
is a monomorphism, for any field $L$. The $A$-homomorphism
$$
{X}\mapr{\wt{\eta}}{}\bigoplus_{i=1}^{m_Z}{E_i}
$$
is a pure monomorphism and we  denote its cokernel by $\wt{Z}$. It
follows that $\L{Z}$ is the cokernel of the map $(*)$ and thus it
is the terminus of the Auslander-Reiten sequence starting at
$\L{\wt{X}}$, for any field $L$.

Then $\wt{\CS}_n\cup\{\wt{Z}\}$ is a set of $\mathbb{Z}$-frames of
$\CS_{n+1}=\CS_n\cup\{Z\}$. \epv

\section{$\mathbb{Z}$-schemes of submodules}\label{sec:grassmann}
Let $A=\ZZ Q/I$ be an order and $N_1,N_2,M$ be $A$-lattices. We
prove that the locally closed sets
$\CE(\K{N_1},\K{N_2};\K{M})\subseteq Gr_{\scriptsize
{\bdim\,N_1}}(\bdim\,M,K)$ are defined "almost" independently on
the base field $K$.

If $F_1,...,F_r$ are homogeneous polynomials in $n+1$ variables
with coefficients in a field $K$, then we denote by $
V_K(F_1,...,F_r)$ the Zariski-closed subset of $\mathbb{P}^{n}(K)$
defined by $F_1=...=F_r=0$.

We denote by ${\frak K}$ the algebraic closure of the field of
rational numbers $\QQ$.

\begin{prop}\punkt Assume that the set $U\subseteq
\mathbb{P}^n({\frak K})$ is invariant under the action of every
automorphism of the field ${\frak K}$. Then the Zariski closure of
$U$ is defined over $\ZZ$, that is,  there are homogeneous
polynomials $F_1,...,F_r \in\ZZ[T_0,...,T_n]$ such that
$$
\ov{U}=V_{{\frak K}}(F_1,...,F_r).
$$
\label{prop:invariant}\end{prop}

{\bf Proof.} Let
$$
\ov{U}=V_{{\frak K}}(H_1,...,H_s),
$$
where $H_1,...,H_s\in{\frak K}[T_0,...,T_n]$ are homogeneous
polynomials. Denote by $L_1$ the subfield of ${\frak K}$ generated
by all coefficients of $H_1,...,H_s$. Let $L$ be the normal
closure of $L_1$ over $\QQ$ and we denote by $G$ the Galois group
of the extension $\QQ\subseteq L$. Assume that
$$
G=\{\sigma_1,...,\sigma_p\},
$$
where $p=[L:\QQ]$.

Given a polynomial $H$ with coefficients in ${\frak K}$ and an
automorphism $\sigma$ of ${\frak K}$ we denote by $H^{\sigma}$ the
polynomial obtained from $H$ by  the action of $\sigma$ on the
coefficients of $H$. Since $U$ is invariant under the
automorphisms of ${\frak K}$, every $H_i^{\sigma_j}$, $i=1,...,s$,
$j=1,...,p$ vanishes on $U$.

Let $s_1,...,s_p$ be the standard symmetric polynomials in $p$
variables. Observe that
$$
H_{i,j}:=s_i(H_j^{\sigma_1},...,H_j^{\sigma_p})
$$
is invariant under each element of $G$, thus its coefficients
belong to $L^G=\QQ$. Clearly, each $H_{i,j}$ vanishes on $U$.
Moreover, since the unique common zero of $s_1,...,s_p$ is
$(0,...,0)$ it follows that
$$
\ov{U}=V_{{\frak K}}(H_{i,j}: i=1,...,s, j=1,...,p).
$$
Therefore $\ov{U}$ is defined by vanishing of polynomials with
coefficients in $\QQ$. Multiplying by denominators we can assume
that
$$
\ov{A}=V_{{\frak K}}(H_1,...,H_q)
$$
for some homogeneous polynomials $H_1,...,H_q$ with coefficients
in $\ZZ$.
 \epv

\begin{cor}\punkt  Let $U\subseteq \mathbb{P}^n({\frak K})$
 be a~locally closed set which is invariant under each
automorphism of the field ${\frak K}$. There are homogeneous
polynomials $F_1,...,F_r, G_1,...,G_s\in\ZZ[T_0,...,T_n]$ such
that
$$
U=V_{{\frak K}}(F_1,...,F_r)\setminus V_{{\frak K}}(G_1,...,G_s).
$$ \label{cor:locclosed}\end{cor}

{\bf Proof.} Corollary follows immediately from Proposition
\ref{prop:invariant} applied to $U$ and $\ov{U}\setminus U$. \epv

Let $B=\mathbb{C}Q/J$ be a~representation directed
$\mathbb{C}$-algebra and assume that $A=\ZZ Q/I$ is the
$\ZZ$-order corresponding to $B$ as in Theorem \ref{thm:Z-form}).

\begin{lem}\punkt Assume that $K$ is a subfield of a field $L$
and $A$ is a $\ZZ$-order.

{\rm (1)}  Let $B$ be a $K$-algebra. If $X,Y$ are $B$-modules then
$X\cong Y$ if and only if the $B\otimes_KL$-modules $X\otimes_KL$
and $Y\otimes_KL$ are isomorphic.

{\rm (2)} If $N_1,N_2,M$ are $A$-lattices then
$$
\CE(\K{N_1},\K{N_2};\K{M})= Gr_{\scriptsize
{\bdim\,N_1}}(\bdim\,M,K) \cap \CE(\L{N_1},\L{N_2};\L{M}).$$
\label{lem:VoverZ}\end{lem}

\textbf{Proof.} The assertion (1) is proved in \cite[Lemma
3.2]{JL1}, whereas (2) follows directly from (1). \epv

\begin{thm}\punkt Let $N_1,N_2,M$ be $A$-lattices.
There exist a finite set $\Sigma$ of prime numbers and there exist
homogeneous polynomials $F_1,...,F_r, G_1,...,G_s$ with integral
coefficients such that
$$
\CE(N_1^{(K)},N_2^{(K)};M^{(K)}) =V_{K}(F_1,...,F_r)\setminus
V_{K}(G_1,...,G_s),
$$
for any field  $K$ of characteristic not belonging to $\Sigma$.
\label{thm:overZ}\end{thm}

 \textbf{Proof.} Since any automorphism of the field ${\frak K}$
is constant on $\mathbb{Z}$ and, by Theorem \ref{thm:Z-form},
$B$-modules are defined over $\mathbb{Z}$, it is easy to see that
the set $\CE(N_1^{({\frak K})},N_2^{({\frak K})};M^{({\frak K})})$
is closed under each automorphism of the field ${\frak K}$.
Therefore it follows from Lemma \ref{lem:locclosed} and Corollary
\ref{cor:locclosed} that there are homogeneous polynomials
$F_1,...,F_r, G_1,...,G_s$such that
$$
\CE(N_1^{({\frak K})},N_2^{({\frak K})};M^{({\frak K})})=V_{{\frak
K}}(F_1,...,F_r)\setminus V_{{\frak K}}(G_1,...,G_s)$$

 Note that the fact $x\in \CE(N_1^{(K)},N_2^{(K)};M^{(K)})$ can be
written as a~first order formula in the language of rings
\cite[Chapter 10]{JL2}. Applying the Characteristic Transfer
Principle \cite[Theorem 1.14]{JL2} to the formula:
\begin{small}$$
x\in \CE(N_1^{(K)},N_2^{(K)};M^{(K)})\leftrightarrow
[(F_1(x)=...=F_r(x)=0)\wedge (G_1(x)\neq 0 \vee ... \vee
G_s(x)\neq 0)].
$$
\end{small}
yields the existence of a finite set $\Sigma$ such that the
equality
$$
\CE(N_1^{({ K})},N_2^{({ K})};M^{({K})})=V_{{
K}}(F_1,...,F_r)\setminus V_{{K}}(G_1,...,G_s)$$ holds whenever
$K$ is an algebraically closed field of characteristic not
belonging to $\Sigma$.

 Thanks to Lemma \ref{lem:VoverZ} the equality holds for every
field $K$ such that $char\; K\notin C$. \epv

\section{Proof of Theorems \ref{thm:main} and \ref{thm:main1}}\label{sec:proof}

In order to prove Theorem \ref{thm:main1} we have to relate the
Euler-Poincar\'{e} characteristic $\chi$  of $\CE(N_1,N_2;M)$ with Hall
polynomials. For the properties of the Euler-Poincar\'{e}
characteristic the reader is referred to \cite{durfee} and
\cite{fulton}. We will need the following lemma, which can
be deduced from  \cite{fulton} (see also \cite{CC}, \cite[Proposition 6.1]{reineke2005} and \cite[Lemma 8.1]{reineke2008}).


\begin{lem}
  Let $F_1,...,F_r,  G_1,...,G_s$ be
homogeneous polynomials with integral coefficients in $p$
variables. Given a field $K$ we denote
$$
V_K=V_K(F_1,...,F_r)\setminus V_K(G_1,...,G_s).
$$
Assume that there exists a polynomial $\phi\in\ZZ[t]$ and a finite
set $C$ of prime numbers such that
$$
|V_K|=\phi(|K|)
$$
for any finite field $K$ such that $char\,K\notin C$.

Then
$$
\chi(V_{\mathbb{C}})=\phi(1).
$$
\epv \label{lem:weil}
\end{lem}

Let $B=\mathbb{C}Q/J$ be a~representation-directed
$\mathbb{C}$-algebra and $A=\mathbb{Z}Q/I$ the $\ZZ$-order
corresponding to $B$ as in  Theorem \ref{thm:Z-form}. The
$K$-algebras $A^{(K)}$ are representation-directed for all fields
$K$ and their Auslander-Reiten quivers coincide. Denote this
common Auslander-Reiten quiver  by $\Gamma$. \vsp

\textbf{Proof of Theorem \ref{thm:main1}.} Let $B=\mathbb{C}Q/J$
be a~representation directed $\mathbb{C}$-algebra and let
$A=\mathbb{Z}Q/I$ be the corresponding $\mathbb{Z}$-order, such
that $A^{(\mathbb{C})}\cong B$. Consider functions
$a,b,c:\Gamma_0\to\mathbb{N}$. By Theorem \ref{thm:hall-poly},
$$|\CE(M(A^{(K)},a),M(A^{(K)},c);M(A^{(K)},b))|=
F_{M(A^{(K)},c),M(A^{(K)},a)}^{M(A^{(K)},b)}=\varphi_{ca}^b(|K|),$$
for any finite field $K$, where $\varphi_{ca}^b$ are Hall
polynomials associated with $a,b,c\in\Gamma$. Thanks to Theorem
\ref{thm:overZ} and Lemma \ref{lem:weil} we obtain
$$\chi(\CE(M(B,a),M(B,c);M(B,b)))=\varphi_{ca}^b(1),$$
and we are done. \epv

\textbf{Proof of Theorem \ref{thm:main}.} Let $B$ be
a~representation-directed $\mathbb{C}$-algebra,
$\ind(B)=\{X_1,\ldots,X_m\}$ and let $\Gamma$ be the
Auslander-Reiten quiver of $B$.
Moreover let $x_1,\ldots,x_m\in \Gamma_0$ be such that $X_i\cong M(B,x_i)$,
for all $i=1,\ldots,m$. Note that the $\mathbb{Z}$-Lie algebras $L(B)$ and $\CK(B)$ are
free $\mathbb{Z}$-modules with basis $v_{X_1},\ldots,v_{X_m}$ and $u_{x_1},\ldots,u_{x_m}$,
respectively. Let $$ F:\CK(B)\to L(B)$$ be the~homomorphism of $\mathbb{Z}$-modules
given by
$$
F(u_{x_i})=\left\{\begin{array}{rl} v_{X_i} &  \mbox{if } \dim_{\mathbb{C}}X_i \mbox{ is odd } \\
     -v_{X_i} & \mbox{otherwise}\end{array} . \right.
$$
We show that $F$ is a~homomorphism of Lie algebras. Let $x_i,x_j\in \Gamma_0$.
Consider the case that there exist $x_t\in\Gamma_0$  and a~short exact sequence
$$ 0\to X_j\to X_t\to X_i\to 0,$$
then, by (\ref{eq:K}) and (\ref{eq:ringel-alg}), 
$$F([u_{x_i},u_{x_j}])=F(\varphi_{x_ix_j}^{x_t}(1)\cdot u_{x_t})=
(-1)^{\dim_{\mathbb{C}}X_t-1}\varphi_{x_ix_j}^{x_t}(1)\cdot v_{X_t}.$$
On the other hand, by (\ref{eq:L}, \ref{eq:Riedtmann-alg} and Theorem \ref{thm:main1}),
$$\begin{array}{lcl}
[F(u_{x_i}),F(u_{x_j})]&=&(-1)^{\dim_{\mathbb{C}}X_i+\dim_{\mathbb{C}}X_i-2}[v_{X_i},v_{X_j}]
\\ &=& 
(-1)^{\dim_{\mathbb{C}}X_i+\dim_{\mathbb{C}}X_i-1}[v_{X_j},v_{X_i}]\\
&=& (-1)^{\dim_{\mathbb{C}}X_t-1}\chi(\CE(X_j,X_i;X_t))\cdot v_{X_t} \\
&=& (-1)^{\dim_{\mathbb{C}}X_t-1}\varphi_{x_ix_j}^{x_t}(1)\cdot v_{X_t}. 
\end{array}
$$
It follows that $F([u_{x_i},u_{x_j}])=[F(u_{x_i}),F(u_{x_j})]$.
In the remaining cases the proof is analogous. Then $F$ is a~homomorphism of Lie algebras.
Obviously $F$ is an~isomorphism of Lie algebras, because it is isomorphism of
free albelian groups $\CK(B)$ and $L(B)$. This
finishes the proof. \epv

\begin{cor}\punkt Let $B$ be a~representation-directed $\mathbb{C}$-algebra.

{\rm (a)} The complex Lie algebras $L(B)\otimes_{\mathbb{Z}}\mathbb{C}$ and
$\CK(B)\otimes_{\mathbb{Z}}\mathbb{C}$ are isomorphic.

{\rm (b)} The algebra $\CH(B)_1\otimes_{\mathbb{Z}}\mathbb{C}$ is the
universal enveloping algebra of $L(B)\otimes_{\mathbb{Z}}\mathbb{C}$. \end{cor}

\textbf{Proof.} The assertion (a) is obvious. By \cite[Proposition
5]{ringel1}, $\CH(B)_1\otimes_{\mathbb{Z}}\mathbb{C}$ is the
universal enveloping algebra of
$\CK(B)\otimes_{\mathbb{Z}}\mathbb{C}$, then  the statement (b)
follows from Theorem \ref{thm:main}. \epv

\begin{rem}\punkt
{\rm Combinatorial properties of the Lie algebra
$\CK(B)\otimes_{\mathbb{Z}}\mathbb{C}$ are investigated in
\cite{kos08}. In particular, the Lie algebra
$\CK(B)\otimes_{\mathbb{Z}}\mathbb{C}$ is described by generators
and relations, for any representation directed
$\mathbb{C}$-algebra $B$.
 This together with Theorem \ref{thm:main} give us a~description of the Lie algebra
$L(B)\otimes_{\mathbb{Z}}\mathbb{C}$ by generators and relations.}
\end{rem}

%


\begin{thebibliography}{999}

\bibitem{ASS} I. Assem,  D. Simson and A. Skowro\'{n}ski, {\em "Elements of
Representation Theory of Associative Algebras", Vol. I: Techniques
of Representation Theory}, London Mathematical Society Student
Texts, 65. Cambridge University Press, Cambridge, 2006.

\bibitem{ARS} M. Auslander, I. Reiten and S. Smal\o, {\em "Representation
theory of Artin algebras"}, Cambridge Studies in Advanced
Mathematics 36,  Cambridge University Press,
  1995.

  \bibitem{Bo1} K. Bongartz,  {\it Zykellose Algebren sind nicht zugellos},
  in Representation theory, II,
     Lecture Notes in Mathematics 832, Springer, Berlin, 1980, 97--102.

\bibitem{bo83}
K. Bongartz, {\it Algebras and quadratic forms}, { J. London Math.
Soc.,} 28 (1983), 461-469.

\bibitem{CC} P. Caldero and F. Chapoton, {\it Cluster algebras as Hall algebras of quiver
representations}, Commentarii Mathematici Helvetici, 81 (2006), 595-616.

\bibitem{DP} P. Dr\"axler and J. A. de la Pe\~na, {\it On the existence of postprojective
components in the Auslander-Reiten quiver of an algebra}, {\rm Tsukuba J.
Math.}  20 (1996), 457-469.

\bibitem{durfee} A. H. Durfee, {\it Algebraic varieties which are a~disjoint union
of subvarieties}, Lecture Notes in Pure and Appl. Math. 105 (1987), 99-102.

\bibitem{fulton}
W. Fulton, {\it "Introduction to Toric Varieties"}, Princeton University Press, 1993.

\bibitem{JL1}  Ch. Jensen and H. Lenzing,  Homological dimension and
representation type of algebras under base field extension,
 {\em Manuscripta Math.}, {\bf 39},  1-13 (1982)

\bibitem{JL2}
C. U. Jensen, H. Lenzing, {\it "Model theoretic algebra with particular emphasis on fields, rings, modules"},
Algebra Logic Appl. 2, Gordon \& Breach, 1989.


\bibitem{Ka} S. Kasjan, Representation-directed algebras form an open scheme,
 {\em Colloq. Math.} {\bf 93} (2002), 237-250.

\bibitem{KaPe} S. Kasjan and J.A. de la Pe\~na, Constructing the
preprojective components of an algebra, {\em J. Algebra} {
179}(1996), 793-807.

\bibitem{kos08} J. Kosakowska, {\it Lie algebras associated with quadratic forms and their applications to Ringel-Hall algebras}, preprint 2008.

\bibitem{reineke2001} M. Reineke, {\it Generic extensions and multiplicative bases of quantum groups at
$q=0$}, An Electronic Journal of the Amer. Math. Soc., Vol. 5
(2001), 147-163.

\bibitem{reineke2005} M. Reineke, {\it Counting rational points of quiver
moduli}, arXiv:math/0 50 5389 v1.

\bibitem{reineke2008} M. Reineke, {\it Moduli of representations of quivers},
arXiv:0802.2147v1 [math.RT].

\bibitem{riedtmann} Ch. Riedtmann, {\it Lie algebras generated by
indecomposables}, J. Algebra 170 (1994), 526-546.

\bibitem{Ri1099} C. M. Ringel, {\it Tame Algebras and Integral Quadratic Forms},
Lecture Notes in Mathematics,  1099 (Springer-Verlag, Berlin,
Heidelbegr, New York, Tokyo 1984).

\bibitem{ringel} C. M. Ringel, {\it Hall algebras and quantum
groups}, Invent. Math. 101 (1990), 583-592.

\bibitem{ringel1} C. M. Ringel, {\it Hall algebras}, Banach Center
Publications, Vol. 26, Warsaw 1990, 433-447.

\bibitem{ringel90} C. M. Ringel, {\it Hall polynomials for the representation-finite hereditary algebras},
Adv. in Math. 84 (1990), 137-178.

\bibitem{ringel89} C. M. Ringel, {\it From representations of quivers via Hall and Loewy algebras to quantum groups}, Proceedings Novosibirsk Conference 1989. Contemporary Mathematics 131.2 (1992), 381-401.


\bibitem{SkWe} A. Skowro\'nski and M. Wenderlich,
Artin algebras with directing indecomposable projective modules.
{\em J. Algebra} {\bf 165} (1994), 507--530.


\end{thebibliography}
\end{document}